\documentclass[letterpaper,12pt]{article}
\usepackage{amssymb}
\usepackage{amsmath}
\usepackage{makeidx}
\usepackage{amscd}

\newtheorem{theo}{Theorem}[section]
\newtheorem{prop}[theo]{Proposition}
\newtheorem{lem}[theo]{Lemma}
\newtheorem{cor}[theo]{Corollary}
\newtheorem{defi}[theo]{Definition}

\def \bfx{{\bf x}}

\def \Br {{\rm{Br}}}

\def \Ga {{\Gamma}}

\def \Pic {{\rm {Pic}}}
\def \Div {{\rm {Div}}}

\def \Gal {{\rm{Gal}}}

\def \AA{{\mathbf A}}
\def \A{{\mathbb A}}
\def \P{{\mathbb P}}
\def \Spec {{\rm{Spec}}}
\def \dim {{\rm{dim}}}
\def \Hom {{\rm {Hom}}}

\def \Pic {{\rm {Pic}}}

\def\ov{\overline}

\def\k{\kappa}
\def \Z {{\mathbb Z}}
\def \Q {{\mathbb Q}}

\def \Ext {{\rm Ext}}

\def \AA {{\rm A}}

\def\G{{\mathbb G}}

\def\calT{{\cal T}}
\def\T{{\cal T}}

\def\lra{\longrightarrow}

\def\O{{\cal O}}

\def\O{{\cal O}}

\def\Ga{\Gamma}

\def\m{{\mathfrak m}}

\newcommand{\bthe}{\begin{theo}}
\newcommand{\ble}{\begin{lem}}
\newcommand{\bpr}{\begin{prop}}
\newcommand{\bco}{\begin{cor}}
\newcommand{\bde}{\begin{defi}}
\newcommand{\ethe}{\end{theo}}
\newcommand{\ele}{\end{lem}}
\newcommand{\epr}{\end{prop}}
\newcommand{\eco}{\end{cor}}
\newcommand{\ede}{\end{defi}}

\title{Descent on toric fibrations}
\author{Alexei N. Skorobogatov}
\date{}
\begin{document}
\baselineskip=15pt
\maketitle

\begin{abstract}
\noindent We describe descent on families of torsors of a constant torus.
A recent result of Browning and Matthiesen then implies
that the Brauer--Manin obstruction controls the Hasse principle and 
weak approximation when the ground field is $\Q$ and the singular 
fibres are all defined over $\Q$.
\end{abstract}


\section{Introduction}

Let $T$ be a torus over a number field $k$.
Let $X$ be a smooth, proper, 
geometrically integral variety with a surjective morphism 
$f:X\to\P^1_k$ whose generic fibre $X_{k(t)}$ is 
geometrically integral and is
birationally equivalent to a $k(t)$-torsor of $T$.
The main result of this note says that the set $X(k)$ is dense in
$X(\AA_k)^\Br$ if certain auxiliary varieties satisfy
the Hasse principle and weak approximation. 

These varieties are given by explicit equations.
Choose $\A^1_k\subset\P^1_k$ so
that the fibre of $f$ at $\infty=\P^1_k\setminus\A^1_k$ is smooth.
Let $P_1,\ldots,P_r$ be closed points of $\P^1_k$ such that each fibre of
the restriction of $f$ to $\P^1_k\setminus (P_1\cup\ldots\cup P_r)$ 
is {\em split}, i.e. contains
an irreducible component of multiplicity 1 that is geometrically
irreducible. By a well known result (Lemma \ref{le1})
the fibre $X_{P_i}$, for $i=1,\ldots, r$, has an
irreducible component of multiplicity $1$.
We fix such a component in each $X_{P_i}$ and define $k_i$
as the algebraic closure of the residue field $k(P_i)$
in the function field of this component. For $i=1,\ldots, r$
let $p_i(t)\in k[t]$ be the monic irreducible polynomial such that
$P_i$ is the zero set of $p_i(t)$ in $\A^1_k$, and let $a_i$
be the image of $t$ in $k(P_i)=k[t]/(p_i(t))$.
Let $u,v$ be independent variables, and let $z_i$ be a $k_i$-variable,
for $i=1,\ldots, r$.
For $\alpha=\{\alpha_i\}$, where $\alpha_i\in k(P_i)^*$, we
define the quasi-affine variety 
$W_\alpha\subset\A^2_k\times \prod_{i=1}^r R_{k_i/k}(\A^1_{k_i})$
by 
\begin{equation}
\alpha_i(u-a_iv)=N_{k_i/k(P_i)}(z_i),\quad (u,v)\not=(0,0), \quad 
i=1,\ldots, r, \label{eqW}
\end{equation}
The varieties $W_\alpha$ are smooth and geometrically irreducible.
Over an algebraic closure of $k$ such a variety is given by
$\sum_{i=1}^r [k(P_i):k]$ equations in 
$2+\sum_{i=1}^r [k_i:k]$ variables.
We can now state our main result.

\bthe \label{one}
Suppose that for each $\alpha \in\prod_{i=1}^r k(P_i)^*$
the variety $W_\alpha$ satisfies the Hasse principle and 
weak approximation. Then $X(k)$ is dense in $X(\AA_k)^\Br$.
\ethe

The results of this kind are obtained by the descent method
of Colliot-Th\'el\`ene and Sansuc, and have a long history.
When the relative dimension and 
the number of singular geometric fibres of $f$
are small, geometric proofs of the Hasse principle
and weak approximation for $W_\alpha$ were obtained by 
Colliot-Th\'el\`ene, Sansuc,
Swinnerton-Dyer and others, see \cite{CSS, CTSal, S90, SD, CTSk} and
\cite[Ch. 7]{Sk}. The analytic tool in these proofs is
Dirichet's theorem on primes in an arithmetic progression.
Over $k=\Q$ one can do more if one uses analytic methods: the
circle method \cite{HS, CHS}, sieve methods \cite{BH, DSW} and 
recent powerful results from additive combinatorics 
\cite{GT1, GT2, GTZ, BMS, HSW, Smeets}. Note that
the circle method can sometimes be applied 
over arbitrary number fields, see \cite{SJ, SchS}. 

When $k=k(P_1)=\ldots=k(P_r)=\Q$, a recent 
result of Browning and Matthiesen obtained by methods of
additive combinatorics \cite[Thm 1.3]{BM}
establishes the Hasse principle and weak approximation for $W_\alpha$.
Hence we deduce the following corollary of Theorem \ref{one}.

\bco \label{main}
Let $X$ be a smooth, proper, geometrically integral variety over $\Q$,
and let $f:X\to\P^1_\Q$ be a surjective morphism satisfying 
the following properties.

{\rm (a)} There is a torus $T$ over $\Q$
such that the generic fibre $X_{\Q(t)}$ of $f$ is birationally equivalent
to a $\Q(t)$-torsor of $T\times_\Q\Q(t)$.

{\rm (b)} There exists a finite subset $E\subset \P^1_\Q(\Q)$
such that $X\setminus f^{-1}(E)\to 
\P^1_\Q\setminus E$ has split fibres.

Then $X(\Q)$ is dense in $X(\AA_\Q)^\Br$.
\eco

This generalises a recent result due to A. Smeets, namely
the unconditional counterpart of \cite[Thm. 1.1, Rem. 1.3]{Smeets}.
For a higher-dimensional version of this statement see Proposition
\ref{higher}.

For a number field $k$ of degree $n=[k:\Q]$ we write $\bfx=(x_1,\ldots,x_n)$
and denote by $N_{k/\Q}(\bfx)$ the norm form
${\rm Norm}_{k/\Q}(x_1\omega_1+\ldots+x_n\omega_n)$, where
$\omega_1,\ldots,\omega_n$ is a basis of $k$ as a vector space
over $\Q$. The following corollary extends \cite[Thm. 1.1]{BM} 
to the case of a product of norm forms. It generalises the unconditional
version of \cite[Cor. 1.6]{Smeets}, a number of statements from 
\cite[Section 4]{HSW} and the main result of
\cite{SchS} in the case of the ground field $\Q$.

\bco \label{co1}
Let $k_1,\ldots,k_n$ be number fields and let
$m_1,\ldots,m_n$ be positive integers with
$gcd(m_1,\ldots,m_n)=1$. Let $L_i\in\Q[t_1,\ldots,t_s]$ 
be polynomials of degree $1$, for $i=1,\ldots,r$.
Let $X$ be a smooth
and proper variety over $\Q$ that is birationally equivalent to the
affine hypersurface
\begin{equation}
\prod_{i=1}^r L_i(t_1,\ldots,t_s)=\prod_{j=1}^n N_{k_j/\Q}(\bfx_j)^{m_j}.
\label{norms}
\end{equation}
Then $X(\Q)$ is dense in $X(\AA_\Q)^\Br$.
\eco

Note that repetitions among $L_1,\ldots,L_r$ are allowed here. 
Corollary \ref{co1}
is a particular case of Proposition \ref{sys} which deals with
several equations like (\ref{norms}) and extends \cite[Thm. 1.3]{BM}.

This note consists of two sections.
In \S2 we make preliminary remarks, some of which, like
Corollary \ref{c}, are not needed in the proof of our main results
but could possibly be of independent interest. In \S3
we prove Theorem \ref{one}, Corollary \ref{co1} and their
generalisations.

Our proof of Theorem \ref{one} uses descent like \cite{BMS},
\cite{DSW} or \cite{SchS}, and not the fibration method
like \cite{Smeets} or \cite{HSW}. It was inspired by the 
approach of Colliot-Th\'el\`ene and Sansuc \cite{CTS}
to degeneration of torsors of tori, and by
their computation of universal torsors 
on conic bundles \cite[Section 2.6]{CS87}. We apply
open descent based on Harari's formal lemma
as in \cite{CTSk}, with an improvement found in \cite{CHS}.

The author is grateful to Jean-Louis Colliot-Th\'el\`ene for 
useful discussions over the past many years. I would like to thank Daniel
Loughran for his question that led me to Proposition \ref{0},
and Olivier Wittenberg for useful comments.

\section{Reduction of varieties 
defined over a dicretely valued field} \label{S2}

Let $R$ be a discrete valuation ring with the field of fractions
$K$, the maximal ideal $\m$ and the residue field $\k=R/\m$. 
We assume that $\k$ is perfect.
Let $j:\Spec(K)\to \Spec(R)$ and $i:\Spec(\k)\to \Spec(R)$ 
be the embeddings of
the generic and the special points, respectively.
We have an exact sequence of \'etale sheaves on $\Spec(R)$:
\begin{equation}
0\to\G_{m,R}\to j_*\G_{m,K} \to i_*\Z_{\k}\to 0,
\label{alef}
\end{equation}
see \cite[Examples II.3.9]{EC}. 
Let $T$ be a torus over $R$. Applying $\Hom_R(\widehat T,\cdot)$
we obtain an exact sequence of abelian groups
\begin{equation}
H^1(R,T)\to H^1(K,T)\to \Ext^1_{\k} (\widehat T,\Z),
\label{bet}
\end{equation}
see \cite[p. 70]{Sk} for a detailed proof. We note that
$\Ext^1_{\k} (\widehat T,\Z)=H^1(\k,\widehat T^\circ)$, where
$T^\circ$ is the $\k$-torus dual to $T$, that is, such that
$\widehat T^\circ=\Hom(\widehat T,\Z)$ as Galois modules.

For an $R$-scheme $X$ we denote by $X_K$ and $X_\k$
the generic and special fibres of $X$, respectively.

\ble \label{1}
Let $Y$ be a $K$-torsor of $T_K$. 
Suppose that there is an integral normal scheme $X$ and
a surjective morphism $X\to\Spec(R)$ with  
integral fibres such that the generic fibre 
$X_K$ is birationally equivalent to $Y$. 

{\rm (i)}  Let $\k'$ be the algebraic closure of $\k$
in the function field of the special fibre $\k(X_\k)$. 
Then the image of $[Y]$ in $H^1(\k,\widehat T^\circ)$
is in the kernel of the restriction map 
$H^1(\k,\widehat T^\circ)\to H^1(\k',\widehat T^\circ)$.

{\rm (ii)} If $X_\k$ is geometrically integral, 
then there is an $R$-torsor $Z$ of $T$
such that $Y\cong Z_K$.
\ele
{\em Proof}. Let $R'$ be the local ring of the special fibre $X_\k$.
Since $X$ is normal and $X_\k$ is integral, 
$R'$ is a discrete valuation ring.
Since $X$ is integral, the field of fractions of $R'$ is $K(X_K)$.
A local parameter of $R$ is also a local parameter of $R'$,
because $X_\k$ is integral. Thus $\m'=\m\otimes_R R'$
is the maximal ideal of $R'$, and the residue field $R'/\m'$
is the field of functions on the special fibre $X_\k$.

The pullback from $R$ to $R'$ gives rise to the commutative diagram
$$\begin{array}{ccccc}
H^1(R',T)&\to &H^1(K(X_K),T)&\to &H^1(\k(X_\k),\widehat T^\circ)\\
\uparrow&&\uparrow&&\uparrow\\
H^1(R,T)&\to &H^1(K,T)&\to& H^1(\k,\widehat T^\circ)
\end{array}
$$
The restriction of the diagonal $Y\to Y\times_K Y$ to the generic fibre 
of $Y$ is a $K(Y)$-point of $Y$.
Hence the torsor $Y$ is split by the field extension $K(X_K)=K(Y)$, so that
the class $[Y]\in H^1(K,T)$ goes to zero in $H^1(K(X_K),T)$.

The right vertical map in the diagram factorises as follows:
$$H^1(\k,\widehat T^\circ)\lra H^1(\k',\widehat T^\circ)
\tilde\lra H^1(\k(X_\k),\widehat T^\circ).$$
The second map here is an isomorphism because 
$\widehat T^\circ$ is a finitely generated free abelian group and $\k'$ 
is algebraically closed in $\k(X_\k)$. This implies (i).

If the $\k$-scheme $X_\k$ is geometrically integral, the field $\k$ is
algebraically closed in $\k(X_\k)$, that is, $\k'=\k$. Then
$[Y]$ is in the image of the map $H^1(R,T)\to H^1(K,T)$, 
and this proves (ii). $\Box$

\medskip

The statement of Lemma \ref{1} (i) leaves open the question
to what extent $\k'$ is determined by the field $K(X_K)$. 
We treat this as a question about integral,
regular, proper schemes over a discretely valued field, see
Corollary \ref{c} below.

For an integral variety $V$ over a field $k$, we write $k_V$ for
the algebraic closure of $k$ in $k(V)$.

\bpr \label{0}
Let $Y$ and $Y'$ be integral regular schemes that are proper over $R$.
If there is a dominant rational map from $Y_K$ to $Y'_K$, then for any 
irreducible component $C\subset Y_\k$ of multiplicity $1$ there exists
an irreducible component $C'\subset Y'_\k$ of multiplicity $1$ such
that $\k_{C'}\subset \k_{C}$.
\epr
{\em Proof.} Write $F=K(Y_K)$ and $F'=K(Y'_K)$. We are given an inclusion
$F'\subset F$, or, equivalently, a morphism $\Spec(F)\to\Spec(F')$.
Let $\O_{C}$ be the local ring of the generic point
of $C$. Since $Y$ is regular, $\O_{C}$ is a discrete valuation
ring. The multiplicity of $C$ is $1$, that is,
the maximal ideal of $\O_{C}$ is $\m\O_{C}=\m\otimes_R\O_{C}$. 
Since $Y$ is integral, the field of fractions of $\O_{C}$ is $F$. 
The residue field of $\O_{C}$ is $\k(C)$.

After the base change from $R$ to $\O_{C}$ we obtain a morphism
$Y'\times_R \O_{C}\to \O_{C}$.
Its generic fibre $Y'_{F}=Y'_K\times_K F$ has an $F$-point
$\Spec(F)\to\Spec(F')\to Y'_{F}$
coming from the $F'$-point defined by the diagonal morphism
$Y'_K\to Y'_K\times_K Y'_K$. 
The morphism $Y'\times_R \O_{C}\to \O_{C}$ is proper, 
and by the valuative criterion
of properness any $F$-point of its generic fibre extends to a section
of the morphism.

Since $C$ has multiplicity $1$,
the special fibre of $Y'\times_R \O_{C}\to \O_{C}$
is $Y'_\k\times_\k \k(C)$. A section of 
$Y'\times_R \O_{C}\to \O_{C}$ thus gives rise to
a $\k(C)$-point of $Y'_\k\times_\k \k(C)$,
which can be viewed also as a $\k(C)$-point of $Y'_\k$.
Since $Y'\times_R \O_{C}$ is regular, 
any section of $Y'\times_R \O_{C}\to \O_{C}$ meets the special fibre 
at a smooth point. In particular, this point 
belongs to exactly one geometric irreducible component,
moreover, this component has multiplicity $1$ (see the calculation
in the proof of \cite[Lemma 1.1 (b)]{S96}, or \cite[Lemme 3.8]{W}). 

Let $U\subset Y'_\k$ be the smooth locus of $Y'_\k$. 
Since $U$ contains a $\k(C)$-point, it is non-empty.
By Stein factorisation, the structure morphism $U\to
\Spec(\k)$ factors through the surjective morphism $U\to
\Spec(L)$ with geometrically connected fibres, where $L$ is 
an \'etale $\k$-algebra. Explicitly, $L$ is
the direct sum of finite field extensions of $\k$, 
each of which is the algebraic closure of $\k$ in the function
field of an irreducible component of $Y'_\k$ of multiplicity $1$.
The image of the composed map $\Spec(\k(C))\to \Spec(L)$ is connected,
hence this is $\Spec(\k_{C'})$ for some irreducible component 
$C'\subset Y'_\k$ of multiplicity $1$. Thus $\k_{C'}\subset \k(C)$
and hence $\k_{C'}\subset \k_{C}$. $\Box$

\bco \label{c}
Let $X$ be an integral regular scheme that is proper over $R$. 
Let $\Sigma_X$ be the partially
ordered set of irreducible components of multiplicity $1$ of $X_\k$,
where $C$ dominates $D$ if $\k_{D}\subset\k_{C}$.
The set of finite field extensions $\k\subset \k_C$, 
where $C$ is a minimal element of $\Sigma_X$, is a birational
invariant of the generic fibre $X_K$.
\eco
{\em Proof.} Suppose that proper $R$-schemes $X$ and $Y$ are integral
and regular with birationally equivalent generic fibres,
that is, $K(X_K)\cong K(Y_K)$. Let $C$ be a minimal element
of $\Sigma_{X}$. By Proposition \ref{0} there exists
$C'\in \Sigma_{Y}$ such that $\k_{C'}\subset \k_{C}$. By the same
proposition, there is $C''\in\Sigma_X$ such that $\k_{C''}\subset \k_{C'}$.
By minimality of $C$ we have $\k_{C''}= \k_C$, hence $\k_C=\k_{C'}$. 
If $C'$ is not minimal in $\Sigma_{Y}$, then, by Proposition \ref{0},
$C$ is not minimal in $\Sigma_X$. $\Box$

\medskip

This set of finite extensions of the residue field 
can be explicitly determined when
the generic fibre is a conic, and, more generally, a Severi--Brauer variety,
or a quadric of dimension 2.

\medskip

\noindent{\bf Remark}. Olivier Wittenberg suggested a somewhat different
approach to Proposition \ref{0} and 
Corollary \ref{c}. Consider a discrete valuation $v:
K(X_K)^*\to\Z$ such that the restriction of $v$ to $K^*$
is the given discrete valuation of $K$. 
Let $R'$ be the valuation ring of $v$, 
and let $\k'$ be the algebraic closure
of $\k$ in the residue field of $R'$. Let us call 
$\Theta_X$ the resulting set
of finite field extensions of $\k$
partially ordered by inclusion. It is clear that $\Theta_X$
is a birational invariant of $X_K$. The discrete valuation associated
to an irreducible component of $X_\k$ of multiplicity $1$ is an
example of such a valuation, so we have an
inclusion of partially ordered sets $\Sigma_X\subset\Theta_X$. 
Since $X$ is proper over $R$ and regular, it can be shown that
any extension of the given discrete valuation of $K$ to
$K(X_K)$ gives rise to a morphism of $R$-schemes 
${\rm Spec}(R')\to X$ that factors through the smooth locus of $X/R$.
Hence $\Sigma_X$ and $\Theta_X$ have the same set of minimal elements,
which is thus a birational invariant of $X_K$.

\section{Torsors over toric fibrations} \label{S3}

Let $k$ be a field of characteristic zero.
Let $\bar k$ be an algebraic closure of $k$, and let
$\Ga_k=\Gal(\bar k/k)$.
For a variety $X$ over $k$ we write $\ov X=X\times_k\bar k$.

Let $T$ be a $k$-torus. We write $\widehat T$ for the $\Ga_k$-module
of characters of $T$.

Let $X$ be a smooth, proper, 
geometrically integral variety with a surjective morphism 
$f:X\to\P^1_k$ and the geometrically integral generic fibre $X_{k(t)}$  
which is birationally equivalent to a $k(t)$-torsor of $T$.

\ble \label{le1}
Each fibre of $X\to\P^1_k$ has an irreducible component
of multiplicity~$1$.
\ele
{\em Proof.} The generic fibre of $\ov X\to\P^1_{\bar k}$ is
birationally equivalent to
a torsor of $\ov T\simeq\G_{m,\bar k}^d$, where $d=\dim(T)$.
By Hilbert's theorem 90 we have $H^1(\bar k(t), T)=
H^1(\bar k(t),\G_{m})^d=0$, so this torsor has
a $\bar k(t)$-point. By the lemma of Lang and Nishimura,
$X_{k(t)}$ has a $\bar k(t)$-point too. By the valuative criterion of 
properness, this point extends to a section of the 
proper morphism $\ov X\to\P^1_{\bar k}$. Since $X$ is smooth,
by a standard argument (see 
the proof of \cite[Lemma 1.1 (b)]{S96}, or \cite{W}) any section intersects
each fibre of $\ov X\to\P^1_{\bar k}$ in an irreducible component 
of multiplicity $1$. The lemma follows. $\Box$

\medskip

\noindent{\em Proof of Theorem} \ref{one}. We keep the notation
of \S1. By Lemma \ref{le1}
there is an irreducible component of multiplicity $1$ 
in each $X_{P_i}$, for $i=1,\ldots, r$. We mark these components.
We also mark a geometrically integral irreducible component 
of multiplicity 1 in each of 
the remaining fibres of $f$. Define $Y\subset X$ as the complement 
to the union of all the unmarked irreducible 
components of the fibres of $f:X\to\P^1_k$. It is clear that $Y$
is a dense open subset of $X$. The restriction of $f$ to 
$Y$ is a surjective morphism $f:Y\to\P^1_k$ with integral fibres, and
with proper and geometrically integral generic fibre $Y_{k(t)}=X_{k(t)}$.
It follows that $\bar k[Y]^*=\bar k^*$ and $\Pic(\ov Y)$ is torsion-free.

Let $\Pic(\ov Y)\to\Pic(Y_{\bar k(t)})$ be the homomorphism
of $\Ga_k$-modules induced by the inclusion of the generic
fibre $Y_{\bar k(t)}$ into $\ov Y$. This homomorphism is surjective
since $Y$ is smooth. Let $S$ be a $k$-torus
defined by the exact sequence of $\Ga_k$-modules
\begin{equation}
0\to \widehat S\to \Pic(\ov Y)\to\Pic(Y_{\bar k(t)})\to 0. \label{odin}
\end{equation}
Thus the abelian group $\widehat S$ is generated by the geometric
irreducible components of the fibres $Y_{P_1},\ldots,Y_{P_r}$.
Recall that a {\em vertical torsor} $\T\to Y$ is a torsor
of $S$ whose type is the injective map $\widehat S\to\Pic(\ov Y)$
from (\ref{odin}). According to \cite[Prop. 4.4.1]{Sk} such torsors
exist. For any vertical torsor we have $\bar k[\calT]^*=\bar k^*$
and the abelian group
$$\Pic(\ov\calT)\cong \Pic(Y_{\bar k(t)})\cong \Pic(X_{\bar k(t)})$$
is torsion-free. It follows that $\Br_1(\calT)/\Br_0(\calT)$ is finite.

Recall that for $i=1,\ldots,r$ we write $k(P_i)$ 
for the residue field of $P_i$ and $k_i$ for the algebraic closure
of $k(P_i)$ in the function field $k(Y_{P_i})$. 
The variety $W_\alpha$ is defined by (\ref{eqW}),
which is also equation (4.34) of \cite{Sk}. 

\ble \label{le2}
The variety $W_\alpha$ is smooth and geometrically irreducible,
and the morphism $\pi:W_\alpha\to\P^1_k$ given by $(u:v)$ is faithfully flat.
We have $\bar k[W_\alpha]^*=\bar k^*$ 
and $\Pic(\ov W_\alpha)=0$. 
\ele
{\em Proof.} The first statement is \cite[Lemma 4.4.5]{Sk}.
The second statement is a straightforward adaptation of
\cite[Lemme 2.6.1]{CS87}. $\Box$

\bpr \label{2.1}
Any vertical torsor over $Y$ is birationally equivalent to 
$W_\alpha\times_k Z$, where $\alpha \in\prod_{i=1}^r k(P_i)^*$ and 
$Z$ is a $k$-torsor of $T$.
\epr
{\em Proof.} The local description of torsors
due to Colliot-Th\'el\`ene and Sansuc, see \cite[Thm. 4.3.1, Cor. 4.4.6]{Sk},
can be stated as follows. 
Let $$U'=\A^1_k\setminus\{P_1,\ldots,P_r\}, \quad U=Y\cap f^{-1}(U'), 
\quad V_\alpha=\pi^{-1}(U').$$ 
Then for any vertical torsor $\T/Y$ there exists an 
$\alpha \in\prod_{i=1}^r k(P_i)^*$ such that
the restriction $\T_U=\T\times_Y U$ is isomorphic to the fibred product
$U\times_{U'}V_\alpha$. Let us write $W=W_\alpha$, $V=V_\alpha$.

Let $j':\Spec(k(t))\to U'$ be the embedding
of the generic point of $\A^1_k$ into the open set $U'$.
Let $i_P:P\to U'$ be the embedding of a closed point $P\subset U'$.
Since $\A^1_k$ is smooth,
there is the following exact sequence of \'etale sheaves on $U'$:
$$
0\to\G_{m,U'}\to j'_*\G_{m,k(t)} \to \oplus_{P\in U'} i_{P*}\Z_{k(P)}
\to 0,
$$
see \cite[Examples II.3.9, III.2.22]{EC}. Similarly to the proof of
Lemma \ref{1}, an application of $\Hom_{U'}(\widehat T,\cdot)$
to this sequence produces an exact sequence
\begin{equation}
H^1(U',T)\to H^1(k(t),T)\to \oplus_{P\in U'} H^1(k(P),\widehat T^\circ).
\label{dva}
\end{equation}
Let $\xi\in H^1(k(t),T)$ be the class of a $k(t)$-torsor of $T$
birationally equivalent to $X_{k(t)}$. The map 
$H^1(k(t),T)\to H^1(k(P),\widehat T^\circ)$ in (\ref{dva}) can be computed
in the local ring $R=O_P$. The fibres of $f:U\to U'$ are geometrically
integral, hence Lemma \ref{1} (ii) implies that this map is trivial. Thus
we see from (\ref{dva}) that $\xi$ comes from some $\xi'\in H^1(U',T)$.
By Lemma \ref{1} (i) applied to $f:Y\to\P^1_k$
the image of $\xi'$ in $H^1(k(P_i),\widehat T^\circ)$ goes to zero in 
$H^1(k_i,\widehat T^\circ)$.

The fibre $W_i=\pi^{-1}(P_i)$ is the product of 
a principal homogeneous space of a $k(P_i)$-torus and the
affine $k(P_i)$-variety  defined by $N_{k_i/k(P_i)}(z_i)=0$. Thus
$W_i$ is integral over $k(P_i)$ and the field
$k_i$ is the algebraic closure of $k(P_i)$ 
in the function field of $W_i$. If we set $P_0=\infty$, then $W_0=\pi^{-1}(P_0)$
is geometrically integral by construction, so the algebraic closure of 
$k(P_0)=k$ in $k(W_0)$ is $k_0=k$.
We have a commutative diagram similar to the
one in the proof of Lemma \ref{1}:
$$\begin{array}{ccccc}
H^1(W,T)&\to &H^1(V,T)&\to &\oplus_{i=0}^r H^1(k(W_i),\widehat T^\circ)\\
\uparrow&&\uparrow&&\uparrow\\
H^1(k,T)&\to &H^1(U'_0,T)&\to& \oplus_{i=0}^r H^1(k(P_i),\widehat T^\circ)
\end{array}$$
By the structure of degenerate fibres of $\pi:W\to\P^1_k$
described above, the image of $\xi'$ in 
$H^1(k(W_i),\widehat T^\circ)$ is zero for $i=0,\ldots,r$.
Now the upper sequence in the diagram shows that 
$\pi^*(\xi')\in H^1(V,T)$ comes from $H^1(W,T)$.

By Lemma \ref{le2} we have $\bar k[W_\alpha]^*=\bar k^*$
and $\Pic(\ov W_\alpha)=0$, and thus the
fundamental exact sequence of Colliot-Th\'el\`ene and Sansuc
(see \cite[Thm. 1.5.1]{CS87} or \cite[Cor. 2.3.9]{Sk})
shows that the natural map $H^1(k,T)\to H^1(W,T)$ is an isomorphism.
It follows that any $W$-torsor of $T$ is the product of $W$ and
a $k$-torsor of $T$. $\Box$

\bigskip

\noindent{\em End of proof of Theorem} \ref{one}. 
Let $(M_v)\in X(\AA_k)^\Br$. 
By a theorem of Grothendieck, $\Br_1(X)$ is naturally a subgroup of $\Br_1(Y)$.
We have $\bar k[Y]^*=\bar k^*$, and this implies that 
$\Br_1(Y)/\Br_0(Y)$ is a subgroup of $H^1(k,\Pic(\ov Y))$,
which is finite because $\Pic(\ov Y)$ is torsion-free. Thus we can use
\cite[Prop.~1.1]{CTSk} (a consequence of Harari's formal lemma) which
says that the natural injective map of topological spaces
$Y(\AA_k)^{\Br_1(Y)}\rightarrow X(\AA_k)^{\Br_1(X)}$
has a dense image. Thus we can assume without loss of generality
that $(M_v)\in Y(\AA_k)^{\Br_1(Y)}$.

The main theorem of the descent theory
of Colliot-Th\'el\`ene and Sansuc states that every point in
$Y(\AA_k)^{\Br_1(Y)}$ is in the image of the map
$\calT_0(\AA_k)\to Y(\AA_k)$, where $\calT_0\to Y$ is a universal 
torsor (see \cite[Section 3]{CS87} or \cite[Thm.~6.1.2(a)]{Sk}).
Thus we can find a point $(N_v)\in\calT_0(\AA_k)$ such that
the image of $N_v$ in $Y$ is $M_v$ for all $v$. 

The structure group of $\calT_0\to Y$ is the N\'eron--Severi torus $T_0$
defined by the property $\widehat T_0=\Pic(\ov Y)$.
The dual of the injective map in (\ref{odin}) is a surjective
morphism of tori $T_0\to S$. Let $T_1$ be its kernel.
Then $\calT=\calT_0/T_1$ is a $Y$-torsor with the structure 
group $S$ whose type is the natural map
$\widehat S\to \Pic(\ov Y)$,
so $\calT$ is a vertical torsor. Since $\calT_0$ is a universal torsor,
we have $\bar k[\calT_0]^*=\bar k^*$ and $\Pic(\ov \calT_0)=0$,
hence $\Br_1(\calT_0)=\Br_0(\calT_0)$. 
Let $(P_v)\in \calT(\AA_k)$ be the image of $(N_v)$. By the functoriality
of the Brauer--Manin pairing we see that 
$(P_v)\in \calT(\AA_k)^{\Br_1(\calT)}$. 

By Proposition \ref{2.1}
the variety $\calT$ is birationally equivalent to $E\times W_\alpha$
for some $\alpha \in\prod_{i=1}^r k(P_i)^*$. 
Let $E^c$ be a smooth compactification of $E$.
Then we have a rational map $g$ from the
smooth variety $\calT$ to the proper variety $E^c$, and 
a rational map $h$ from $\calT$ to $W_\alpha$.
By the valuative criterion of properness there is an
open subset $\Omega\subset\calT$ with complement 
of codimension at least 2 in $\calT$ such that $g$ is a 
morphism $\Omega\to E^c$. 
Let $\Omega'\subset\Omega$ be a dense open subset
such that $(g,h)$ defines an open embedding $\Omega'\subset E^c\times W_\alpha$.

By Grothendieck's purity theorem the natural
restriction map $\Br(\calT)\to\Br(\Omega)$ is an isomorphism.
Thus $g^*\Br_1(E^c)\subset\Br_1(\Omega)$, so that 
$g^*\Br_1(E^c)\subset\Br_1(\calT)$. (This argument is borrowed
from \cite{CTSk}.)

Since $\Br_1(\calT)/\Br_0(\calT)$ is finite,
by a small deformation we can assume $P_v\in \Omega(k_v)$ 
for each $v$. Thus $(g(P_v))$ is a well defined element of
$E^c(\AA_k)^{\Br_1(E^c)}$. By Sansuc's theorem,
$E(k)$ is then non-empty and, moreover, is dense in $E^c(\AA_k)^{\Br_1(E^c)}$.

Let $\Sigma$ be a finite set of places
of $k$ containing all the places where we need to approximate.
By assumption we can find a $k$-point in $\Omega'$
that is arbitrarily close to $(g(P_v),h(P_v))$ for $v\in \Sigma$.
We conclude that there is a $k$-point in $\calT$ that is arbitrarily
close to $P_v$ for $v\in \Sigma$. The image of this point in $Y$ approximates
$(M_v)$. This finishes the proof of Theorem \ref{one}. $\Box$

\bpr \label{higher}
Let $X$ be a smooth, proper, geometrically integral variety over $\Q$,
and let $f:X\to\P^n_\Q$ be a surjective morphism satisfying 
the following properties.

{\rm (a)} There is a torus $T$ over $\Q$
such that the generic fibre $X_{\Q(\eta)}$ of $f$ is birationally equivalent
to a $\Q(\eta)$-torsor of $T\times_\Q\Q(\eta)$.

{\rm (b)} There exist hyperplanes 
$H_1,\ldots,H_r\subset \P^n_\Q$ such that $f$ has split fibres 
at all points of codimension $1$ of $\P^n_\Q$ other than $H_1,\ldots,H_r$.

Then $X(\Q)$ is dense in $X(\AA_\Q)^\Br$.
\epr
{\em Proof.} We deduce this from a fibration
theorem due to D. Harari \cite[Thm. 3.2.1]{harari}. Choose a point
$M\in \P^n_\Q(\Q)$ such that the fibre $X_M$ is smooth.
The projective lines through $M$ are in a natural
bijection with the points of $\P^{n-1}_\Q$. Thus an appropriate
open subset $V\subset X$ is a quasi-projective variety equipped 
with a surjective morphism $p:V\to \A^{n-1}_\Q$ with split fibres.
Any point of $X_M(\ov\Q)$ defines a section of $p$ over $\ov\Q$.
The generic fibre of $p$ is birationally equivalent to
a family of torsors of $T$ over an open subset of the projective line,
hence it is geometrically integral and geometrically
rational. In particular, a smooth and proper model of the generic fibre has
trivial geometric Brauer group and torsion-free geometric Picard group.
The fibres of $p$ over $\Q$-points of a non-empty open subset of 
$\A^{n-1}_\Q$ satisfy the assumptions 
of Corollary \ref{main}, hence all the conditions of 
\cite[Thm. 3.2.1]{harari} are satisfied. 
An application of this result proves the proposition.
$\Box$

\bpr \label{sys}
Let $k_1,\ldots,k_n$ be number fields, and let
$L_i\in\Q[t_1,\ldots,t_s]$ be polynomials of degree $1$, for $i=1,\ldots,r$.
Let $m_{ij}\geq 0$, for $i=1,\ldots,\ell$ and $j=1,\ldots, n$,
be integers such that the sublattice of $\Z^n$
generated by the rows of the matrix $(m_{ij})$ is primitive.
Finally, let $d_{hi}\geq 0$ be integers, where $h=1,\ldots,\ell$
and $i=1,\ldots,r$. 
Let $X$ be a smooth
and proper variety over $\Q$ that is birationally equivalent to the
affine variety given by the system of equations
\begin{equation}
c_h\prod_{i=1}^r L_i(t_1,\ldots,t_s)^{d_{hi}}=\prod_{j=1}^n N_{k_j/\Q}(\bfx_j)^{m_{hj}},
\quad h=1,\ldots,\ell, \label{sysnorms}
\end{equation}
where $c_h\in\Q^*$. Then $X(\Q)$ is dense in $X(\AA_\Q)^\Br$.
\epr
{\em Proof.} The condition on the matrix $(m_{ij})$ implies that
the affine variety given by 
$$\prod_{j=1}^n N_{k_j/\Q}(\bfx_j)^{m_{hj}}=1, \quad h=1,\ldots,\ell, $$
is a torus. Let us call it $T$. 
The affine variety $Y$ given by
(\ref{sysnorms}) has a morphism $g:Y\to\A^s_\Q$ given by the coordinates
$t_1,\ldots,t_s$. Let $H_i$ be the hyperplane $L_i=0$. The
restriction of $g$ to $\A^s_\Q\setminus (H_1\cup\ldots\cup H_r)$
is a torsor of $T$. We can choose
a smooth compactification $X$ of the smooth locus $Y_{\rm sm}$
in such a way that there is a surjective morphism $f:X\to\P^s_\Q$
extending $g:Y_{\rm sm}\to\A^s_\Q$. Now we apply
Proposition \ref{higher}. $\Box$

\medskip

For $\ell=1$ this statement is Corollary \ref{co1}.

\bigskip

\noindent Department of Mathematics, South Kensington Campus,
Imperial College London, SW7 2BZ England, U.K. - and - 
Institute for the Information Transmission Problems,
Russian Academy of Sciences, 19 Bolshoi Karetnyi, Moscow, 127994
Russia
\medskip

\noindent{\tt a.skorobogatov@imperial.ac.uk}

\end{document}